\RequirePackage{fix-cm}
\documentclass[smallextended]{svjour3}       
\smartqed  
\usepackage{amsfonts,color,graphicx}
\begin{document}

\title{Robust iteration methods for complex systems with an indefinite matrix term 
}
\subtitle{}


\author{Owe Axelsson, Maeddeh Pourbagher, Davod Khojasteh Salkuyeh}


\institute{
	O. Axelsson \at
	The Czech Academy of Sciences, Institute of Geonics, Ostrava, Czech Republic\\
    \email{owe.axelsson@it.uu.se} \\
    \\
    M. Pourbagher \at
                  Faculty of Mathematical Sciences, University of Guilan,  Rasht, Iran\\
                  \email{mapourbagher@gmail.com}   \\
\\	
	D. K. Salkuyeh \at
              Faculty of Mathematical Sciences, University of Guilan, Rasht, Iran\\
              Center of Excellence for Mathematical Modelling, Optimization and Combinational
              Computing (MMOCC), University of Guilan, Rasht, Iran\\
              \email{khojasteh@guilan.ac.ir}           
}

\date{Received: date / Accepted: date}

\maketitle

\begin{abstract}
Complex valued systems with an indefinite matrix term arise in important applications such as for certain time-harmonic partial differential equations such as the Maxwell's equation and for the Helmholtz equation.
Complex systems with symmetric positive definite matrices can be solved readily by rewriting the complex matrix system in two-by-two block matrix form with real matrices which can be efficiently solved by iteration using the preconditioned square block (PRESB) preconditioning method and preferably accelerated by the Chebyshev method. The appearances of an indefinite matrix term causes however some difficulties. To handle this we propose different forms of matrix splitting methods, with or without any parameters involved. A matrix spectral analyses is presented followed by extensive numerical comparisons of various forms of the methods.   
	
\end{abstract}

\section{Introduction}
\label{sec01}
Complex matrix systems arise in time-harmonic PDE problems such as the linear equation $$\frac{\partial \widetilde{u}}{\partial t} + \mathcal{L} \widetilde{u} = \widetilde{f},$$ defined in a spacial bounded domain $\Omega$ and time interval $[0,T]$, with given boundary conditions \cite{Yousept12}. Here the solution and the given source function have time-harmonic forms
\[
\widetilde{u}=ue^{\frac{i\omega t}{T}},\quad \widetilde{f}=f e^{\frac{i\omega t}{T}},
\] 
where the frequency $\omega$ is a multiple of $2 \pi$. This leads to a complex system of linear equations
\[
i \omega u + \mathcal{L} u = f.
\]
After application of a finite element or finite difference method to the above equation, the following system of linear equations is obtained 
\begin{equation}\label{eq1}
(i \omega M + K) u = f,
\end{equation}
where $M$ and $K$ are two square matrices, typically inertia and stiffness matrices, respectively. When the matrices $M$ and $K$ are symmetric positive semidefinite with at least one of them being positive definite, there are several iteration methods for solving the system \cite{Axelsson-Kucherov,Axelsson201859,MHSS,PMHSS,Salkuyeh-Med,Salkuyeh-MC,Salkuyeh-NLWA,Salkuyeh-NA,Salkuyeh-arXiv,Salkuyeh-IJCM,Salkuyeh-Calcolo,TSS-DKS}. Some of them are directly applied to the main system (see \cite{Salkuyeh-arXiv,Salkuyeh-IJCM,Salkuyeh-Calcolo,TSS-DKS,MHSS,PMHSS}) and some to its real form \cite{Axelsson-Kucherov,Axelsson201859,Salkuyeh-Med,Salkuyeh-MC,Salkuyeh-NLWA}.     

To avoid complex arithmetics one can  rewrite the system (\ref{eq1}) in two-by-two block matrix form with real matrices and vectors \cite{Axelsson-Kucherov,AxLu19} 
$$
\mathcal{A} = \left(\begin{array}{cc}
K & - \omega M \cr \omega M & K
\end{array}\right) \left(\begin{array}{cc}
u_1 \cr u_2
\end{array}\right) = \left(\begin{array}{cc}
f_1 \cr f_2
\end{array}\right).
$$ 	
This can be solved efficiently by iteration using the PRESB preconditioning matrix as follows
$$
\mathcal{B} = \left(\begin{array}{cc}
K & - \omega M \cr \omega M & K + 2 \omega M
\end{array}\right) =\left( \begin{array}{cc}
I & -I \cr 0 & I
\end{array}\right)  
\left(\begin{array}{cc}
K+ \omega M & 0 \cr \omega M & K + \omega M
\end{array}\right)
\left(\begin{array}{cc}
I & I \cr 0 & I
\end{array}\right),
$$
where it is seen that each action of the inverse of the preconditioning matrix just involves two solutions with the real and positive definite matrix $K+\omega M$ (see \cite{AxNeyStroem,ANA17,Axelsson201859,PourbagherJJIAM}). Due to that the Krylov subspace methods normally delete the larger eigenvalues first, such methods exhibit a superlinearly rate of convergence \cite{AKM10}. Further, an elementary spectral analysis shows that the eigenvalues of $\mathcal{B}^{-1}\mathcal{A}$ are contained in the interval $[ \frac{1}{2},1]$.

As shown in \cite{AxLiang}, here it is best to use the Chebyshev semi-iteration as acceleration method. This avoids rounding errors as would arise in the vector orthogonalizations if a Krylov subspace method is used. Furthermore, due to the tight eigenvalue bounds the Chebyshev semi-iteration method converges equally fast as the Krylov method. As has been shown in \cite{LiAxZh20,ANA17}, it is much more reliable when solving ill-posed inner systems, which is illustrated clearly by a test example.

Consider a system with a complex matrix 
\begin{equation}\label{eq2}
Ax:=(W+iT)x = b,
\end{equation}
where $W=W_1-W_2$, with both $W_1$ and $W_2$ being symmetric positive definite, as well as $T$. Hence $W$ is indefinite in general. Li and Wu in \cite{Wu2017313} proposed the modified positive/negative stable splitting (MPNS) method when $W$ is symmetric indefinite. Recently, Pourbagher and Salkuyeh in \cite{PourbagherJJIAM} have presented the symmetric positive definite and negative stable splitting (SNSS) method for solving (\ref{eq2}).
This method can be written as follows
\begin{equation}\label{snss}
\left\{ \begin{array}{l}
(\alpha T+W_{2})                  x^{(k+\frac{1} {2})}  =   ((\alpha +i) T+ W_{1})   x^{(k)  }  -  b, \\
(i(\beta +1)  T  +  W_{1}   )   x^{(k+1)}                 =    (i \beta T + W_{2} )      x^{(k+\frac{1}{2})}  +   b,
\end{array}\right.
\end{equation}
where $x^{(0)}$ is an initial guess and   $ \alpha,\beta>0$ are two parameters to be chosen. It was shown the method is convergent under suitable conditions. The SNSS method induces the preconditioner
$$B_{\alpha,\beta}=\frac{1}{(\alpha- i \beta)} (\alpha T + W_{2}) T^{-1} ( i (\beta +1) T+W_{1}) ,$$
for the system (\ref{Ex2}). In each iteration of the GMRES for solving the system (\ref{eq2}) with the SNSS preconditioner, two subsystems with coefficient matrices  $\alpha T + W_{2}$ and $i (\beta +1) T+W_{1}$ should be solved, which can done using the GMRES method or the Chebyshev semi-iteration method incorporated with the PRESB preconditioner (see \cite{PourbagherJJIAM,Axelsson201859} for more details). Numerical results presented in \cite{PourbagherJJIAM}  show that the SNSS method outperforms the MPNS method.

If the matrix $W$ is indefinite the PRESB preconditioner can not be directly applied to the real form of the system (\ref{eq2}), but one can use some form of matrix splittings, similar to the classical alternating iteration method, see e.g. \cite{Axbook}. Here the PRESB method can be used for each of the block matrices arising in the matrix splitting.   
	We present various forms of such methods with spectral analysis. The methods involve coupled inner-outer iterations. It is shown that the number of the outer iterations are insensitive to the inner solution accuracy, which can be utilized to lower the solution cost.

Throughout the paper we use the following notation. $\|.\|$ denotes the Euclidean norm. Spectral radius of a square matrix is denoted by $\rho(.)$. The imaginary unit is shown by $i$~ ($i=\sqrt{-1}$).  The Kronecker product is denoted by $\otimes$.

This paper is organized as following. In Section \ref{sec02}, we present  three matrix splitting methods for solving the system (\ref{eq2}) and investigate their convergence properties. Some numerical experiments are given in  Section \ref{sec03}. Section \ref{sec04} is devoted to some concluding remarks.

\section{Matrix splitting methods}\label{sec02}


As remarked in the introduction and shown earlier e.g. in \cite{AxNeyStroem,AFN16_2}, there exists an efficient preconditioning method, the PRESB method when solving complex matrix systems with spd matrices rewritten in a two-by-two real block matrix form. 
We show here how this method can be applied also for complex systems with an indefinite matrix. Use then a block iteration method with the splittings, $(W_1 + iT)$, $W_2$ respectively $W_2 - iT$, $W_1$ in the iterative solution method.

Hence, given an initial approximation $x^{(0)}$, for $k=0,1,\ldots$ until convergence, solve
\begin{equation}\label{eq3}
\left\{ \begin{array}{rcl}
(W_1 + iT)x^{(k+\frac{1}{2})} &=& W_2 x^{(k)} + b \\
(W_2 - iT)x^{(k+1)} &=& W_1 x^{(k+\frac{1}{2})} - b.
\end{array}
\right.	
\end{equation}   
We refer to this method as {\bf Method I}. Each iteration involves solving two complex valued subsystems. Both of them can be efficiently solved using the GMRES method or the Chebyshev  semi-iteration method in conjunction with PRESB preconditioner. We shall see that the rate of convergence is insensitive to the accuracy used for these inner systems.

Computing $x^{(k+\frac{1}{2})}$ from the first equation in (\ref{eq3}) and substituting in the second one gives the stationary iteration
\[
x^{(k+1)} = Bx^{(k)}-g , \quad  k=0,1,\ldots,
\]
where 	
\begin{eqnarray*}
	B &=& (W_2 - iT)^{-1} W_1 (W_1 + iT)^{-1} W_2, \\
	g &=& i(W_2-iT)^{-1} T (W_1+iT)^{-1}b.
\end{eqnarray*}
It is easy to verify that if 
\begin{eqnarray*}
	M &=& -\frac{1}{i} (W_1 + iT)T^{-1} (W_2 - iT), \\
	N &=& -\frac{1}{i} W_1 T^{-1} W_2,
\end{eqnarray*}
then $A=M-N$ and $B=M^{-1}N$. 

\begin{theorem}\label{Thm1}
	The iteration (\ref{eq3}) converges unconditionally to the exact solution of (\ref{eq2}) for any initial guess $x^{(0)}$. 
	\begin{proof}
		The iteration matrix of Method I is similar to
		\begin{eqnarray*}
			W_2 B W_2^{-1} &=& W_2 (W_2 - iT)^{-1}W_1(W_1 + iT)^{-1} \\ 
			&=& (I - i T W_2^{-1})^{-1}(I + i T W_1^{-1})^{-1}.
		\end{eqnarray*}
		Here we make a similarity transformation with $T^{1/2}$ to get 
		$$
		\widehat{B} = T^{1/2} W_2 B W_2^{-1}T^{-1/2} = (I-i \widehat{W}_2^{-1})^{-1}(I + i \widehat{W}_1^{-1})^{-1},
		$$
		where $\widehat{W}_i = T^{-1/2}W_i T^{-1/2}$, $i=1,2$. It follows that	
		the absolute values of the eigenvalues of the matrix $B$ are bounded as
		\begin{equation}\label{eq4}
		|\lambda(B)| \leq \|\widehat{B}\| \leq \sqrt{(1+\|\widehat{W}_2\|^{-2})^{-1}(1+\|\widehat{W}_1\|^{-2})^{-1}}<1,
		\end{equation}
		where we have used $\|\widehat{B}\|=\sqrt{\rho(B^* B)}$, and the fact that
		$$
		\|I-i\widehat{W}_i^{-1}\| = \sqrt{1+\|\widehat{W}_i\|^{-2}}, \qquad i=1,2.
		$$
		Hence, convergence of Method I follows from  (\ref{eq4}). 
	\end{proof}	
\end{theorem}


Eq. (\ref{eq4}) leads to an acceptable rate of convergence unless both of $\|\widehat{W}_1\|$ and $\|\widehat{W}_2\|$ are very large. This means that $T$ should be comparable in size to at least one of the matrices $W_1$ or $W_2$. If $\|\widehat{W}_i\|\leq 1$ for at least one of them, then $\|B\| \leq 1/\sqrt{2}$.

From Theorem \ref{Thm1} we deduce that the eigenvalues of the matrix $M^{-1}A$ are clustered in 
a circle with radius 1, centered at $(1,0)$. This means that the Krylov subspace methods like GMRES are suitable for solving the preconditioned system $M^{-1}Ax=M^{-1}b$. In each iteration of GMRES with the  preconditioner $M$  two subsystems with the coefficient matrices $W_1 + iT$ and $W_2 - iT$ should be solved, which can be accomplished using the GMRES method or the Chebyshev  semi-iteration method in conjunction with PRESB preconditioner. 

There is an alternative form of the matrix splitting method. We rewrite (\ref{eq1}) in the form   
$$
(T-i(W_1-W_2))x = -ib,
$$
where $T$ and $W_1$, $i=1,2$ are spd. Here the splitting
\begin{equation}\label{eq5}
\left\{ \begin{array}{rcl}
(T-iW_1)x^{(k+\frac{1}{2})} &=& -iW_2 x^{(k)} - ib \\
(T + iW_2)x^{(k+1)} &=& iW_1 x^{(k+\frac{1}{2})}-ib.
\end{array}
\right.	
\end{equation}
can be used. The corresponding iteration matrix equals
\begin{equation}\label{eq6}
B = (T+iW_2)^{-1} W_1 (T-iW_1)^{-1}W_2.
\end{equation}
Using the similarity transformations as before it follows that this leads to the same bound as in (\ref{eq4}). We refer to this method as {\bf Method II}. Similar to Method I, we can see that the Method II serves the preconditioner
\[
M=-\frac{1}{i}(T-iW_1) T^{-1} (T+iW_2) 
\] 
for the system (\ref{eq2}).

To improve this method we use the following scaling in the matrix splitting method,
\begin{equation}\label{MetIII}
\left\{ \begin{array}{rcl}
(\alpha T+iW_2)x^{(k+\frac{1}{2})} &=& ((\alpha-1)T + i W_1) x^{(k)} - ib \\
(\alpha T - iW_1)x^{(k+1)} &=& ((\alpha-1)T - iW_2) x^{(k+\frac{1}{2})} - ib.
\end{array}
\right.	
\end{equation}
where the scaling parameter $\alpha \geq 1$. We refer to this method as {\bf Method III}. Note that for $\alpha =1$ we get the same form as in (\ref{eq5}). 

Similar to Method I the iteration (\ref{MetIII}) can be written 
\[
x^{(k+1)} = Bx^{(k)}-g , \quad  k=0,1,\ldots,
\]
where 	
\begin{eqnarray*}
	B &=&  (\alpha T - iW_1)^{-1} \left((\alpha - 1)T-iW_2\right)(\alpha T + i W_2)^{-1} \left((\alpha-1)T+iW_1\right), \\
	g &=& (2\alpha-1)i(\alpha T+iW_2)T^{-1}(\alpha T-iW_1)b.
\end{eqnarray*}
It is straightforward to prove that if we define 
\begin{eqnarray*}
	M &=& \frac{1}{(1-2\alpha)i} (\alpha T+iW_2)T^{-1} (\alpha T - iW_1), \\
	N &=& \frac{1}{(1-2\alpha)i} ((1-\alpha) T-iW_2)T^{-1} ((1-\alpha) T + iW_1),
\end{eqnarray*}
then $A=M-N$ and $B=M^{-1}N$. So this method induces the preconditioner $M$ for the system (\ref{eq2}). 

Using a similarity transformation as above we find that 
$$
\|B\| \leq \sqrt{\frac{(\alpha-1)^2 + \|\widehat{W}_2\|^{2}}{ \alpha^2 + \|\widehat{W}_2\|^{2}}}\ \sqrt{\frac{(\alpha-1)^2 + \|\widehat{W}_1\|^{2}}{ \alpha^2 + \|\widehat{W}_1\|^{2}}}.
$$
If for instance $\|\widehat{W}_1\| > \|\widehat{W}_2\|$, we use the upper bound, 
$$
\|B\|^2 \leq \frac{(\alpha -1)^2 + \|\widehat{W}_2\|^2}{\alpha^2 + \|\widehat{W}_2\|^2} = 1- \frac{2\alpha -1}{\alpha^2 + \|\widehat{W}_2\|^2}.
$$
It follows that this upper bound is minimized for $\alpha = \frac{1}{2} + \sqrt{\frac{1}{2} + \|W_2\|^2}$. Hence for large values of $\|\widehat{W}_2\|$, we shall choose $\alpha \approx \|\widehat{W}_2\|$ and then 
$$
\|B\| \approx \sqrt{1 - \frac{\alpha}{\alpha^2 + \|\widehat{W}_2\|^{2}}} \approx 1 - \frac{1}{4} \|\widehat{W}_2\|^{-1}.
$$
We note that the upper bound taken for (\ref{eq6}), i.e. corresponding to $\alpha = 1$, equals 
$$
\|B\| \leq 1/\sqrt{1+\|\widehat{W}_2\|^{-2}} \approx 1- \frac{1}{2}\|\widehat{W}_2\|^{-2}.
$$ 
Hence by choosing $\alpha$ close to $\|\widehat{W}_2\|$, the rate of convergence when $\|\widehat{W}_2 \|$ is large is improved by an order of magnitude. 

As we shall see, for one of the examples, Example 2, we can take $\alpha = \frac{\sigma_1}{\sigma_2}$ if $\sigma_1 \geq \sigma_2$. For example, if $\sigma_1/\sigma_2 = 100$, we can expect a magnitude $10^{-4}$ of iterations for $\alpha = 1$, but only $O(10^{-2})$ for the scaled version. However, due to influence of other factors, this is not always seen. Note here that the number of iterations is approximately  $\frac{1}{\delta}\ln 1/\varepsilon$ if $\|B\| = 1-\delta$ and the stopping tolerance equals $\varepsilon \ll 1$.      


It follows that the iteration method is robust. The complex systems are best solved as shown in the Introduction. Using the PRESB method accelerated by the Chebyshev method, each outer iteration requires inner iterations to solve two systems, each with matrix, $W_1 + T$ or $W_2 + T$. For this it is mostly efficient to use an algebraic multigrid preconditioning method, see e.g. \cite{Not10}. 

It follows that there are three levels of iterations, the outer matrix split method, the PRESB and the innermost iterations to solve the basic matrix systems, $T+W_i$, $i=1,2$. Further, it is most efficient to use a flexible form of outer iteration, see \cite{Saad93,AxVas91}. However, in the numerical tests a direct solver is used for the innermost systems.

\section{Numerical experiments}
\label{sec03}

ln this section, we consider the following three examples for our numerical tests. 
\begin{example}\label{Ex1} \rm
	We consider the complex symmetric linear system of equations
	\begin{equation}
	[(-\omega^2 M + K) + i\omega C]x = b,
	\end{equation}
	where $M$ and $K$ are the inertia and stiffness matrices, respectively. We take $C = \omega C_V + C_H$ where $C_V$ and $C_H$ are the viscous and hysteretic damping matrices, respectively; and $\omega$ is the driving circular frequency. For such time-periodic problems a MINRES solver has been used in \cite{KoKo13}. In our numerical experiments, we set $M = I$, $C_V = 5M$ and $C_H = \mu K$ with a damping coefficient $\mu=0.02$ and $K$ the five-point centered difference matrix approximating the negative Laplacian operator with homogeneous Dirichlet boundary conditions, on a uniform mesh in the unit square $[0, 1] \times [0, 1]$ with the mesh size $h = 1/(m + 1)$. ln this case, we have
	$$
	K = (I \otimes V_m + V_m \otimes I) \in \mathbb{R}^{n \times n},
	$$
	with $V_m = h^{-2}{\rm tridiag}(-1,2,-1) \in \mathbb{R}^{m \times m}$. Hence, the total number of variables is $n=m^2$. In addition, the right-hand side vector $f$ is adjusted such that $b=(1+i)Ae$ where $e=(1,1,\ldots,1)^T \in \mathbb{R}^n$ so the exact solution equals $(1+i)e$. Note that we have $W_1 = K$, $W_2 = \omega^2 M$ and $T = \omega C$. Hence $\|T^{-1/2}W_2 T^{-1/2}\| \approx \|M/(5M+ 0.02 K)^{-1}\|$ is bounded with respect to $\omega$ and not large. Therefore for this problem no scaling is needed. 
\end{example}

\begin{example}\label{Ex2} \rm 	  
	We consider the complex Helmholtz equation in 2-D of the form
	\begin{equation} 	  
	\left\{\begin{array}{rll}
	- \Delta u - \sigma_1 u + i \sigma_2 u = f, & {\rm in} & \Omega,\\
	u = g, & {\rm on} & \partial\Omega,
	\end{array}
	\right.	
	\end{equation}	  
	where
	$$
	\Delta = \sum_{j=1}^{2} \frac{\partial^2}{\partial x_j^2},
	$$
	$\sigma_1 \geq 0 $, $\sigma_2 \geq 0$ and $\Omega = [0,1]^2$. The discretization of the equation above in 2-D, using the second order central difference scheme on an $(m + 2) \times (m+2)$ grid of $\Omega$ with mesh-size $h = 1/(m+1)$ leads to a system of linear equations with coefficients matrix $A = W + iT \in \mathbb{C}^{n \times n}$, such that $n=m^2$, and  
	$$
	W = K - \sigma_1 h^2 (I_m \otimes I_m) \qquad {\rm and} \qquad T = \sigma_2 h^2  (I_m \otimes I_m),
	$$
	with $K = I_m \otimes V_m + V_m \otimes I_m$ and $V_m = {\rm tridiag} (-1,2,-1) \in \mathbb{R}^{m \times m}$. This leads to an indefinite matrix $W$. In addition, the right-hand side vector $b$ is adjusted such that $b=(1+i)Ae$ where $e = (1,1,\ldots,1)^T \in \mathbb{R}^n$. Hence, the exact solution equals $(1+i)e$. In our numerical experiments, we choose different combinations of $(\sigma_1,\sigma_2)$. Note that we let $W_1 = K$, $W_2 = \sigma_1 h^2 (I_m \otimes I_m)$. In this example we include a test with scalar parameter $\alpha=\sigma_1/\sigma_2$ if it larger then unity, otherwise we take $\alpha=1$.
\end{example}

\begin{example}\label{Ex3} \rm 	  
	In this example we consider the complex Helmholtz equation in 2-D of the form
	\begin{equation} 	  
	\left\{\begin{array}{rll}
	- \Delta u - 100 u + 10 i u = e^{x+iy}, & {\rm in} & \Omega,\\
	u = 0, ~~~~~& {\rm on} & \partial\Omega,
	\end{array}
	\right.	
	\end{equation}	  
	where $\Omega=[0,1]\times [0,1]$. Similar to the previous example we discretize the equation with $\sigma_1=100$ and $\sigma_2=10$. Here, it is noted that the right-hand side vector is different from the one obtained in the previous example.  
\end{example}

We solve the examples considered using the preconditioned GMRES or flexible GMRES (FGMRES) \cite{Saad93} methods  in conjunction with the preconditioners of Method~I, Method II, Method III and the SNSS method.
In Examples \ref{Ex1} and \ref{Ex3} the main system is solved using the FGMRES method  and in Example \ref{Ex2} by the GMRES method.
The outer iteration (GMRES of FGMRES) is stopped as soon as the Euclidean norm of the residual of  original system is reduced by a factor of $10^{10}$ and the maximum number of iterations is set to be 1000.
We always solve the inner systems using the Chebyshev semi-iteration method with the PRESB preconditioner. 
In Examples \ref{Ex1} and \ref{Ex3} the iteration of Chebyshev semi-iteration method is terminated as soon as the Euclidean norm of the residual of the un-preconditioned system is reduced by a factor of $10^2$ and in Example \ref{Ex2} by a factor of  $10^{10}$. The maximum number of iteration of Chebyshev semi-iteration is set to be 20. We always use a null vector as an initial guess. It is noted that the innermost subsystems are solved exactly using the sparse Cholesky factorization incorporated with the symmetric approximate minimum degree permutation. To do so, the {\sf symamd.m} command of  \textsc{Matlab} is applied.

All of the numerical experiments are performed in \textsc{MATLAB} R2018b by using a Laptop with 2.50 GHz central processing unit (Intel(R) Core(TM) i5-7200U), 6GB RAM and Windows 10.

Numerical results are presented  in Tables \ref{Tbl1}-\ref{Tbl11}. In these tables we report the number of outer iterations (``Iters"), elapsed CPU time in seconds (``CPU") and the following values
$$
R_k = \frac{\|b-Ax^{(k)}\|}{\|b\|}, \qquad E_k = \frac{\|x^*-x^{(k)}\|}{\|x^*\|}
$$    
to demonstrate the accuracy of the computed solutions, where $x^{(k)}$ and $x^*$ are the computed solution at iteration $k$ and the exact solution, respectively. Table \ref{Tbl6} reports the numerical results of Exmaple \ref{Ex3}  where the exact solution is not known.

As the numerical results show all the preconditioners significantly reduce the number of iterations of GMRES or FGMRES. From the CPU time point of view, we observe that the preconditioners  considerably reduce the CPU time of GMRES or FGMRES, except for large values of omega and for not too small values of the mesh size $h$ (see Table \ref{Tbl1}). This is because then an inertia term determines the major eigenvalues of $A$. This holds in particular when the exact solution is smooth.

The numerical results show that, the CPU times increase somewhat faster than with the optimal value $h^{-2}$ when $h$ decreases. This is due to that an incomplete Cholesky method has been used for the innermost matrix systems. We also see that the new preconditioners can compete with the SNSS preconditioner from the CPU time and the number of iterations point of view. Tables \ref{Tbl7} and \ref{Tbl8}  report the numerical results of FGMRES when the inner tolerance varies from $10^{-2}$ to $10^{-10}$. We observe that by use of a less strict inner tolerance, the number of inner iterations decrease saving CPU-times without any increase of the number of outer iterations.
 
 Tables \ref{Tbl9} and \ref{Tbl10} display the numerical results of GMRES-Method III for the inner and outer tolerances $10^{-10}$ for Example 2 and different values of $\sigma_{1}$ and $\sigma_{2}$ when the value of the parameter $alpha$ varies. As we see GMRES-Method III is not sensitive with respect to the parameter $\alpha$. 

In Table \ref{Tbl11} we list average number of the iterations of Chebyshev semi-iteration for solving the inner systems with the PRESB preconditioner versus the inner tolerance for Example \ref{Ex1} and preconditioner I with $\omega=1$ and $m=256$. It is seen that they vary quite regularly with the tolerance criterion.

In general it is seen that the performance of the methods is quite similar. However Methods I-III do not need any a priori estimations of a parameter and the corresponding iteration method converge unconditionally.


\begin{table}[!ht]
	
	\caption{Numerical results of flexible GMRES method for Example 1 with the inner tolerance $10^{-2}$ and the outer tolerance $10^{-10}$. \label{Tbl1}}
	
	\begin{tabular}{|c|c|c|c|c|c|c|c|c|c|} \hline
		$n=16384$	 &$\omega$  & \multicolumn{1}{c|}{$1$} & \multicolumn{1}{c|}{$5$} & \multicolumn{1}{c|}{$10$} &  \multicolumn{1}{c|}{$15$} & \multicolumn{1}{c|}{$20$}  & \multicolumn{1}{c|}{$25$}\\ \hline
		
		&Iter &   266   &   247    &    203     &    158   &  121   &  98    \\
		No-pre& CPU  &  13.33  &  11.57  &   7.92     &   4.89  &  2.95  &   2.03  \\ 
		&$R_{k}$&    8.53e-11    &    9.53e-11   &   9.36e-11   &   9.01e-11   &   8.86e-11   &   6.37e-11  \\ 
		&$E_k$&  1.95e-10    &   3.45e-10   &  1.20e-10  &   7.99e-11   &   5.94e-11    &  2.91e-11   \\ \hline

		&Iter &     7     &  8  &    8   &   7   &   7   &    7   \\
		FGMRES-Method I& CPU  &  0.45    &   0.50    &      0.51    &  0.50    &  0.50   &  0.48    \\ 
		&$R_{k}$&   2.25e-11   &   1.16e-11    &  6.89e-12    &   5.76e-11    &   1.92e-11   &  9.94e-12    \\ 
		&$E_k$&    3.46e-10   &    1.05e-10    &   2.64e-11  &   1.14e-10  &    2.52e-11    &   9.76e-12  \\ \hline

		& Iter &   7    &   8   &   8    &   7     &   7    &  7    \\
		FGMRES-Method II  &  CPU &   0.45    &   0.44    &  0.43    &  0.40    &  0.41     &  0.39     \\ 			
		&  $R_{k}$     &    4.02e-11    &   1.06e-11    &   6.23e-12    &  6.26e-11    &    2.74e-11    &   1.13e-11     \\ 
		& $E_k$&    3.79e-10   &  1.32e-10   &   2.38e-11   &   1.12e-10   &    3.43e-11    &  1.15e-11   \\ \hline

		&$\alpha_{est}$ &  48.8499   &    8.8514   &   3.8534   & 2.1886 &  1.3571   &  0.8589   \\
		&$\beta_{est}$ &  0.0006   &   0.0033   &    0.0065    &  0.0096   &  0.0126   &  0.0156    \\
		FGMRES-SNSS& Iter  &   6   &   8   &   9  &  9   &  9   &  10   \\
		& CPU &   0.45   &  0.45  &  0.51   &   0.54    &   0.54    & 0.56  \\
		&$R_{k}$ &    2.05e-12   &  1.26e-11   &   9.42e-12    & 2.80e-11  &   4.48e-11   &  1.01e-11    \\
		& $E_k$ &   4.22e-11     &   1.36e-10   &   3.24e-11    &   5.19e-11   &  5.39e-11   &  9.28e-12 \\  \hline \hline 
		
		$n=16384$ &$\omega$  & \multicolumn{1}{c|}{$50$} & \multicolumn{1}{c|}{$100$} & \multicolumn{1}{c|}{$150$} &  \multicolumn{1}{c|}{$200$} & \multicolumn{1}{c|}{$250$}  & \multicolumn{1}{c|}{$300$} \\ \hline

		&Iter &    50     &   29    &    22     &    19   &  17   &  15    \\
		No-pre& CPU  &  0.68   &  0.30  &   0.20     &   0.16  &  0.15  &  0.13   \\ 
		&$R_{k}$&   1.15e-10     &    7.53e-11   &   9.85e-11   &  6.45e-11    &   5.25e-11   &   1.01e-10  \\ 
		&$E_k$&    3.21e-11   &   2.70e-11    &  4.23e-11  &   3.16e-11   &    2.82e-11   &   5.84e-11  \\ \hline

		&Iter &     6     &   5  &    5   &   5   &   5   &   5    \\
		FGMRES-Method I& CPU  &   0.45   &    0.40   &    0.39      &   0.40   &  0.38   &  0.39    \\ 
		&$R_{k}$&   1.11e-11   &    3.26e-11   &   3.97e-12   &    2.10e-12   &  1.03e-12    &   6.75e-13   \\ 
		&$E_k$&     9.51e-12   &   2.45e-11     &  3.49e-12   &   1.77e-12  &   6.37e-13     &   4.94e-13  \\ \hline

		& Iter &    6   &   5   &    5   &     5   &    5   &  5    \\
		FGMRES-Method II  &  CPU &     0.34  &   0.30    &   0.31   &   0.30   &   0.29    &  0.30     \\ 			
		&  $R_{k}$     &     1.15e-11    &  6.38e-11     &   1.88e-11    &   1.04e-11   &   6.31e-12    &    3.71e-12    \\ 
		& $E_k$&    9.33e-12   &  4.99e-11   &   1.64e-11   &   9.00e-12   &   5.38e-12     &  3.12e-12   \\ \hline

		&$\alpha_{est}$ &  0.13262   &   0.6181    &  0.7730    & 0.8465 &  0.8883   &  0.9147   \\
		&$\beta_{est}$ &   0.0291  &   0.0513   &    0.0686    &  0.0825   &  0.0938   &  0.1031    \\
		FGMRES-SNSS& Iter  &  14    &  9    &   8  &  8   &  7   &  7   \\
		& CPU &   0.68   &  0.53  &  0.48   &    0.49   &   0.45    &  0.47 \\
		&$R_{k}$ &   6.94e-11    &  1.27e-11   &   2.88e-11    &  7.05e-12 &   4.73e-11   &   2.87e-11   \\
		& $E_k$ &    4.64e-11    &  8.35e-12    &    2.01e-11   &   5.01e-12   &   3.43e-11  & 2.11e-11  \\  \hline
		
	\end{tabular}
	
\end{table}

\begin{table}[!ht]
	
	\caption{Numerical results of flexible GMRES for Example 1 with the inner tolerance $10^{-2}$ and the outer tolerance $10^{-10}$. \label{Tbl2}}
	
	
	\begin{tabular}{|c|c|c|c|c|c|c|c|c|c|} \hline
		$n=65536$	 &$\omega$  & \multicolumn{1}{c|}{$1$} & \multicolumn{1}{c|}{$5$} & \multicolumn{1}{c|}{$10$} &  \multicolumn{1}{c|}{$15$} & \multicolumn{1}{c|}{$20$}  & \multicolumn{1}{c|}{$25$}\\ \hline

		&Iter &     503    &   471    &    398     &   308    &   237  &   186   \\
		No-pre& CPU  &  432.69   & 375.50   &    267.16     &  161.64   &  95.97  &   59.80  \\ 
		&$R_{k}$&    9.94e-11    &   9.31e-11    &   9.47e-11   &   9.49e-11   &   9.12e-11   &  9.90e-11   \\ 
		&$E_k$&   9.86e-10    &   3.88e-10    &  2.09e-10  &   1.51e-10   &   9.56e-11    &   1.46e-10  \\ \hline

		&Iter &       7   &   8  &    8   &    7  &    7  &  7     \\
		FGMRES-Method I& CPU  &   1.79   &    2.10   &     2.20     &    2.16  &  2.17   &  2.19    \\ 
		&$R_{k}$&   6.23e-12   &   6.47e-12    &  5.18e-12    &   2.67e-11    &   1.07e-11   &   5.59e-12   \\ 
		&$E_k$&    2.18e-10   &    1.76e-10    &  5.33e-11   &  1.30e-10   &    3.20e-11    &   1.22e-11  \\ \hline

		& Iter &    7   &   8   &    7   &    7    &    7   &   7   \\
		FGMRES-Method II  &  CPU &  1.90     &   1.88    &  1.76    &  1.66    &  1.68     &   1.70    \\ 			
		&  $R_{k}$     &    1.70e-11     &   4.95e-12    &   6.15e-11    &   3.23e-11   &   1.53e-11    &    6.99e-12    \\ 
		& $E_k$&   4.39e-10    &  1.67e-10   &  6.02e-10    &   1.42e-10   &    4.06e-11    &  1.35e-11   \\ \hline

		&$\alpha_{est}$ &  48.9621   &    8.9622   &   3.9623   & 2.2958 &   1.4626  &  0.9627   \\
		&$\beta_{est}$ &  0.0001   &   0.0008   &   0.0016    &  0.0025   &  0.0033   &  0.0041    \\
		FGMRES-SNSS& Iter  &   6   &  8    &  8   &   9  &  9   &   9  \\
		& CPU &   1.84   &  2.02  &  1.89   &   2.07    &   1.95    &  2.01 \\
		&$R_{k}$ &   2.09e-12    &  6.87e-12   &   8.70e-11    & 1.58e-11  &  3.27e-11    &  8.39e-11    \\
		& $E_k$ &   1.13e-10     &   1.62e-10   &   7.89e-10    &   7.67e-11   &  9.44e-11   &  1.59e-10 \\  \hline
		\hline
		$n=65536$ &$\omega$  & \multicolumn{1}{c|}{$50$} & \multicolumn{1}{c|}{$100$} & \multicolumn{1}{c|}{$150$} &  \multicolumn{1}{c|}{$200$} & \multicolumn{1}{c|}{$250$}  & \multicolumn{1}{c|}{$300$}\\ \hline
		
		&Iter &    101     &    59   &     45    &    38   &  33   &  30    \\
		No-pre& CPU  &  26.27   &  6.55  &    4.05    &  2.93   &  2.27  &   1.88  \\ 
		&$R_{k}$&    9.93e-11    &   7.59e-11    &  8.34e-11    &   7.35e-11  &   9.77e-11   &  8.11e-11  \\ 
		&$E_k$&   3.37e-11    &   2.25e-11    &  2.62e-11  &   2.52e-11   &    3.51e-11   &  3.15e-11   \\ \hline

		&Iter &     6     &   5  &   5    &   5   &   5   &  5     \\
		FGMRES-Method I& CPU  &   2.00   &   1.77    &    1.78      &  1.75    &  1.76   &  1.72    \\ 
		&$R_{k}$&   1.29e-11   &   2.97e-11    &   1.49e-11   &   3.44e-12    &   4.96e-12   &   5.50e-12   \\ 
		&$E_k$&    1.20e-11   &    1.36e-11    &  9.16e-12   &   1.78e-12  &   4.19e-12     &   2.93e-12  \\ \hline

		& Iter &    6   &  5    &   5    &    5    &   5    &  5    \\
		FGMRES-Method II  &  CPU &   1.50    &   1.25    &  1.26    &  1.25    &   1.27    &  1.26     \\ 			
		&  $R_{k}$     &    1.38e-11     &   7.80e-11    &   3.04e-11    &  1.44e-11   &   9.69e-12    &   6.85e-12     \\ 
		& $E_k$&   1.25e-11    &  6.03e-11   &   2.16e-11   &   1.17e-11   &    8.04e-12    &  5.75e-12   \\ \hline

		&$\alpha_{opt}$ &   5  &   5    &   5   & 5 &  5   &  5   \\
		&$\beta_{opt}$ &  0.0081   &   0.0157   &   0.0227     &  0.0293   &  0.0354   &  0.0412    \\
		FGMRES-SNSS& Iter  &   8   &  8    &   8  &  8   &  8   &  7   \\
		& CPU &   1.85   &  1.87  &  1.85   &   1.95   &   1.80    &  1.62 \\
		&$R_{k}$ &   6.47e-11    &  2.60e-11   &   1.33e-11    & 8.34e-12  &  5.53e-12    &   7.67e-11   \\
		& $E_k$ &    5.74e-11    &   1.84e-11   &   8.68e-12    &   5.43e-12   &  3.57e-12   &  4.91e-11 \\  \hline

	\end{tabular}
\end{table}

\begin{table}[!ht]
	
	\caption{Numerical results of GMRES  for Example 2 for $\sigma_{1}=1000$ and $\sigma_{2}=10$ with the inner and outer tolerances $10^{-10}$. \label{Tbl3}}
	
	\begin{tabular}{|c|c|c|c|c|c|c|c|c|c|c|} \hline
		&$n$ & 4096 & 16384  &  65536  \\ \hline
		
		&Iter &  233     &  464    &   933      \\
		No-pre& CPU  &  3.05   &   36.37     &  512.96   \\ 
		&$R_{k}$&  9.29e-11  &  9.55e-11    &  9.94e-11    \\ 
		&$E_k$&  8.35e-11   &   2.96e-10  & 6.57e-10    \\ \hline

		&Iter &    67   &   67     &   67       \\
		GMRES-Method I& CPU  &   2.05    &   7.85     &   32.24   \\ 
		&$R_{k}$&  9.05e-11    &    9.17e-11      &    9.70e-11      \\ 
		&$E_k$&    3.47e-09    &  3.28e-09  &   3.74e-09     \\ \hline
		
		&$\alpha$& 1 & 1  & 1   \\
		&Iter &   66    &   66    &    67      \\
		GMRES-Method III& CPU  &  2.20    &   8.50     &  38.70   \\ 
		&$R_{k}$&  9.20e-11   &    6.56e-11     &    7.08e-11     \\ 
		&$E_k$&   5.66e-09    &  3.61e-09  &   3.89e-09    \\ \hline

		&$\alpha_{est}$& 100 & 100  & 100   \\
		&Iter &   59 &  59    &   60      \\
		GMRES-Method III & CPU  &  1.70   &  6.80    &  31.00  \\ 
		&$R_{k}$& 6.17e-11  &   7.40e-11       &  9.14e-11    \\ 
		&$E_k$&    8.92e-11   &   1.05e-10  &   1.47e-10   \\ \hline

		&$\alpha_{opt}$ &  10  &   5   & 5  \\
		&$\beta_{opt}$ &  1   & 0.9   &   0.9   \\
		GMRES-SNSS& Iter  &  66  & 66  & 67 \\
		& CPU &  1.66  & 7.24  &   30.00  \\
		&$R_{k}$ & 8.76e-11  &  6.81e-11   &  7.96e-11  \\
		& $E_k$ &  4.64e-09   &  3.57e-09  & 3.56e-09  \\  \hline
		
	\end{tabular}
\end{table}

\begin{table}[!ht]
	
	\caption{Numerical results of GMRES for Example 2 for $\sigma_{1}=100$ and $\sigma_{2}=100$ with the inner and outer tolerances $10^{-10}$. \label{Tbl4}}
	\begin{tabular}{|c|c|c|c|c|c|c|c|c|c|c|} \hline
		&$n$ & 4096 & 16384  &  65536  \\ \hline

		&Iter &   139    &  273    &    527     \\
		No-pre& CPU  &   1.30  &   12.83     &   186.47  \\ 
		&$R_{k}$&  7.29e-11  &  9.30e-11    &  9.82e-11    \\ 
		&$E_k$&  3.82e-11   &  1.53e-10   &  1.77e-09   \\ \hline
		
		&Iter &   12    &   12     &    12      \\
		GMRES-Method I& CPU  &   0.50    &   1.55     &  6.59    \\ 
		&$R_{k}$&  1.11e-11    &    1.45e-11      &  1.59e-11       \\ 
		&$E_k$&    2.67e-11    &  3.04e-11  &    3.20e-11    \\ \hline

		&$\alpha =\alpha_{est}$& 1 & 1  & 1   \\
		&Iter &   12    &   12    &   12       \\
		GMRES-Method III& CPU  &  0.48    &   1.70     &  8.00   \\ 
		&$R_{k}$&  3.59e-11  &    1.45e-11     &    3.74e-11     \\ 
		&$E_k$&  3.88e-11     &  3.89e-11  &  4.05e-11     \\ \hline

		&$\alpha_{opt}$ &  5  &  5    & 5  \\
		&$\beta_{opt}$ & 0.1    & 0.1   &  0.1    \\
		GMRES-SNSS& Iter  &  12  & 12  & 12\\
		& CPU & 0.45   &  1.50 &  6.02   \\
		&$R_{k}$ & 1.26e-11  &  1.65e-11   & 1.81e-11   \\
		& $E_k$ &  3.69e-11   &  3.98e-11  & 4.11e-11  \\  \hline

	\end{tabular}
\end{table}

\begin{table}[!ht]	
	\caption{Numerical results of GMRES for Example 2 for $\sigma_{1}=100$ and $\sigma_{2}=10$ with the inner and outer tolerances $10^{-10}$. \label{Tbl5}}
	\begin{tabular}{|c|c|c|c|c|c|c|c|c|c|c|} \hline
		&$n$ & 4096 & 16384  &  65536  \\ \hline

		&Iter &   148    &   291   &   573      \\
		No-pre& CPU  &  1.44   &   14.19     &  212.70   \\ 
		&$R_{k}$&  8.30e-11  &  9.67e-11    &  9.49e-11    \\ 
		&$E_k$&  5.08e-11   &  1.29e-10   &  3.64e-10   \\ \hline
		
		&Iter &    13   &    13    &   13       \\
		GMRES-Method I& CPU  &  0.50     &  1.77      &   7.56   \\ 
		&$R_{k}$&   3.33e-11   &    3.95e-11    &    8.91e-11     \\ 
		&$E_k$&     1.06e-10   &  1.23e-10  &   2.26e-10     \\ \hline
		
		&$\alpha$& 1 & 1  & 1   \\
		&Iter &    13   &    13   &  14        \\
		GMRES-Method III& CPU  &  0.53    &  1.95      & 8.01    \\ 
		&$R_{k}$&  3.37e-11   &   5.52e-11      &  1.58e-11      \\ 
		&$E_k$&   1.00e-10    &  1.63e-10  &  4.66e-11     \\ \hline

		&$\alpha_{est}$& 10 & 10  & 10   \\
		&Iter &  14  &  14    &   15      \\
		GMRES-Method III & CPU  &  0.58   &  2.15    &  9.03  \\ 
		&$R_{k}$& 6.27e-11  &    8.07e-11     &  3.69e-12     \\ 
		&$E_k$&   5.36e-11    &  6.63e-11   &  2.76e-10    \\ \hline

		&$\alpha_{opt}$ &  5  &  5    &  5 \\
		&$\beta_{opt}$ &  0.1   &  0.1  &  0.1    \\
		GMRES-SNSS& Iter  & 13   & 13  & 14\\
		& CPU &  0.53  & 1.71  &   7.14  \\
		&$R_{k}$ & 3.01e-11  &   2.76e-11    &  5.44e-11  \\
		& $E_k$ &  9.90e-11   & 9.22e-11   &  6.38e-11 \\  \hline

	\end{tabular}
\end{table}

\begin{table}[!ht]	
	\caption{Numerical results of flexible GMRES for Example \ref{Ex3}.\label{Tbl6}}
	\begin{tabular}{|c|c|c|c|c|c|c|c|c|c|c|} \hline
		&$n$ & 4096 & 16384  &  65536  \\ \hline

		&Iter    &   290       &   576       &   $\dag$      \\
		No-pre&   CPU  &   5.09      &  62.55      &     -     \\ 
		&$R_{k}$ &  8.63e-11   &  9.51e-11   &     - \\ \hline
		
		& Iter        &    25       &    25          &   25       \\
		FGMRES-Method I    & CPU         &   0.28      &   0.76         &   3.92   \\ 
		& $R_{k}$     &   2.18e-11  &     2.65e-11   &    2.85e-11  \\ \hline
		
		& $\alpha$    &     1       &      1         &  1  \\
		FGMRES-Method III  & Iter        &    24       &      24        & 24   \\
		& CPU         &    0.36    &    1.10        &  5.48       \\
		& $R_{k}$     &  3.28e-11    &  3.88e-11     & 4.10e-11   \\ \hline

		& $\alpha$    &     10      &      10        &  10  \\
		FGMRES-Method III  & Iter        &    27       &      28        & 28   \\
		& CPU         &    0.40     &    1.22       &  6.33        \\
		& $R_{k}$     &  9.98e-11   &  2.55e-11      & 2.91e-11   \\ \hline
		
		& $\alpha_{opt}$ &  5        &       5        &    5 \\
		& $\beta_{opt}$  &  0.1      &       0.1      &  0.1    \\
		FGMRES-SNSS   & Iter           &  25       &       25       &   25 \\
		& CPU            &  0.35     &       0.90     &   3.80  \\
		&$R_{k}$         & 1.76e-11  &   2.02e-11  &  2.12e-11 \\ \hline
		
	\end{tabular}
\end{table}

\begin{table}[!ht]
	
	\caption{Numerical results of FGMRES-Method II for different inner tolerances Example 1 and $n=16384$.} \label{Tbl7}
	\begin{tabular}{|c|c|c|c|c|c|c|c|} 
		\hline
		Inner  & $\omega$ & 1 & 5  & 10 & 15 & 20 & 25  \\ 
		tolerance  & & & & & & & \\ \hline
		
		& Iter &   7    &   8   &   8    &   7     &   7    &  7    \\
		1e-02  &  CPU &   0.45    &   0.44    &  0.43    &  0.40    &  0.41     &  0.39     \\ 			
		&  $R_{k}$     &    4.02e-11    &   1.06e-11    &   6.23e-12    &  6.26e-11    &    2.74e-11    &   1.13e-11     \\ 
		& $E_k$&    3.79e-10   &  1.32e-10   &   2.38e-11   &   1.12e-10   &    3.43e-11    &  1.15e-11   \\ \hline
		
		& Iter &      7     &    8    &    8   &    7    &   7    &   7    \\
		1e-04  &  CPU &   0.58   &   0.65    &   0.64   &  0.55    &   0.57    &  0.57       \\ 			
		&  $R_{k}$     &   3.00e-11     &   1.17e-11    &   5.87e-12    &  5.44e-11  &    2.10e-11    &   7.92e-12     \\ 
		& $E_k$&   2.78e-10    &  1.43e-10   &   2.28e-11   &   1.01e-10   &   2.62e-11    &  7.94e-12   \\ \hline

		& Iter &     7      &     8   &    8   &   7     &   7    &   7    \\
		1e-06  &  CPU &  0.75    &   0.80    &   0.84   &  0.75    &   0.74    &   0.75      \\ 			
		&  $R_{k}$     &    3.00e-11    &   1.17e-11    &   5.89e-12    &  5.46e-11  &   2.10e-11     &    7.94e-12    \\ 
		& $E_k$&   2.78e-10    &   1.43e-10  &   2.28e-11   &   1.02e-10   &   2.62e-11    &  7.96e-12   \\ \hline

		& Iter &       7    &    8    &    8   &    7    &   7    &   7    \\
		1e-08  &  CPU &   0.90   &   1.01    &  1.00    &  0.91    &   0.89    &   0.90      \\ 			
		&  $R_{k}$     &    3.00e-11    &   1.17e-11    &   5.89e-12    &  5.46e-11  &    2.10e-11    &   7.94e-12     \\ 
		& $E_k$&   2.78e-10    &  1.43e-10   &   2.28e-11   &   1.02e-10   &   2.62e-11    &   7.96e-12  \\ \hline

		& Iter &      7     &    8    &   8    &    7    &   7    &   7    \\
		1e-10  &  CPU &   1.05   &   1.15    &   1.18   &  1.04    &   1.03    &   1.05      \\ 			
		&  $R_{k}$     &    3.00e-11    &   1.17e-11    &   5.89e-12    &  5.46e-11  &   2.10e-11     &   7.94e-12     \\ 
		& $E_k$&   2.78e-10    &  1.43e-10  &   2.28e-11   &  1.02e-10    &   2.62e-11    &  7.96e-12   \\ \hline

	\end{tabular}
\end{table}

\begin{table}[!ht]
	
	\caption{Numerical results of FGMRES-Method II for different inner tolerances Example 1 and $n=16384$.} \label{Tbl8}
	\begin{tabular}{|c|c|c|c|c|c|c|c|} 
		\hline
		Inner  & $\omega$ & 50 & 100  & 150 & 200 & 250 & 300  \\ 
		tolerance  & & & & & & & \\ \hline
		
		& Iter &    6   &   5   &    5   &     5   &    5   &  5    \\
		1e-02  &  CPU &     0.34  &   0.30    &   0.31   &   0.30   &   0.29    &  0.30     \\ 			
		&  $R_{k}$     &     1.15e-11    &  6.38e-11     &   1.88e-11    &   1.04e-11   &   6.31e-12    &    3.71e-12    \\ 
		& $E_k$&    9.33e-12   &  4.99e-11   &   1.64e-11   &   9.00e-12   &   5.38e-12     &  3.12e-12   \\ \hline	
		
		& Iter &      6     &    5    &   5    &   4     &   4    &  4     \\
		1e-04  &  CPU &  0.48    &   0.45    &   0.44   &  0.38    &    0.38   &   0.37      \\ 			
		&  $R_{k}$     &   7.09e-12     &   9.61e-12    &   9.19e-13    &  5.50e-11  &    1.96e-11    &  8.08e-12      \\ 
		& $E_k$&   5.57e-12    &   7.71e-12  &   7.55e-13   &   4.60e-11   &   1.64e-11    &  6.80e-12   \\ \hline
		
		& Iter &     6      &   5     &   5    &   4     &   4    &  4     \\
		1e-06  &  CPU &  0.63    &   0.58    &   0.59   &   0.47   &   0.49    &    0.49     \\ 			
		&  $R_{k}$     &    7.12e-12    &   9.64e-12    &   9.24e-13    &  5.51e-11  &    1.95e-11    &   8.04e-12     \\ 
		& $E_k$&    5.59e-12   &  7.73e-12   &   7.58e-13   &   4.61e-11   &    1.64e-11   &  6.78e-12   \\ \hline
		
		& Iter &     6      &   5     &   5    &   4     &   4    &   4    \\
		1e-08  &  CPU &   0.77   &   0.68    &   0.68   &  0.60    &   0.55    &   0.54      \\ 			
		&  $R_{k}$     &   7.12e-12     &   9.64e-12    &   9.24e-13    &  5.51e-11  &   1.95e-11     &    8.04e-12    \\ 
		& $E_k$&   5.59e-12    &   7.73e-12  &   7.58e-13   &   4.61e-11   &   1.64e-11    &  6.78e-12   \\ \hline
		
		& Iter &    6       &    5    &    5   &   4     &    4   &   4    \\
		1e-10  &  CPU &  0.93    &   0.82    &  0.80    &   0.65   &   0.65    &   0.66      \\ 			
		&  $R_{k}$     &    7.12e-12    &   9.64e-12    &   9.24e-13    &  5.51e-11  &    1.95e-11    &   8.04e-12     \\ 
		& $E_k$&   5.59e-12    &  7.73e-12   &   7.58e-13   &   4.61e-11   &   1.64e-11    &   6.78e-12  \\ \hline
		
	\end{tabular}
\end{table}

\begin{table}[!ht]
	
	\caption{Numerical results of GMRES-Method III for the inner and outer tolerances $10^{-10}$ for Example 2, $\sigma_{1}=1000$ and $\sigma_{2}=10$.} \label{Tbl9}
	\label{exact}
	\begin{tabular}{|c|c|c|c|c|c|c|} 
		\hline
		$\alpha$ & 1 & 5 & 10  & 20 & 30 & 40  \\ 
		$n=65536$ & & & & & &\\ \hline
		
		Iter &    67    &   68     &   68   &   68    &   65    &   62    \\
		CPU &   39.68   &    39.06   &   39.83  &   39.73  &   38.70    &  35.12       \\ 			
		$R_{k}$     &   7.08e-11    &  3.88e-11  &   5.60e-11    &  4.12e-11   &  6.77e-11  &  9.68e-11     \\ 
		$E_k$&    3.89e-09   & 7.66e-10   &  7.02e-10   &   3.10e-10  &   3.37e-10   &  3.84e-10  \\ \hline
		
		$\alpha$& 50&  60 & 70 & 80 & 90 & 100 \\
		$n=65536$  & & & & & &\\ \hline
		
		Iter &    63    &    59    &   59   &   59    &    59   &   60    \\
		CPU &    37.53  &   34.08    &  33.63   &  33.61   &   32.43    &   35.04      \\ 			
		$R_{k}$     &   9.68e-11   &  6.54e-11  &   8.34e-11    &  9.38e-11   &  7.44e-11  &   9.14e-11     \\ 
		$E_k$&   1.53e-10    &  1.91e-10  &  2.08e-10   &  1.73e-10  &   1.24e-10   &  1.47e-10  \\ \hline
		
	\end{tabular}
\end{table}

\begin{table}[!ht]
	
	\caption{Numerical results of GMRES-Method III for the inner and outer tolerances $10^{-10}$ for Example 2, $\sigma_{1}=100$ and $\sigma_{2}=10$.} \label{Tbl10}
	\begin{tabular}{|c|c|c|c|c|c|} 
		\hline
		$\alpha$ & 1 & 2 & 3  & 4 & 5   \\ 
		$n=65536$ & & & & & \\ \hline
		
		Iter &     14   &    13   &   13   &    13  &  17   \\
		CPU &  9.39   &   8.74   & 8.99   &  8.85  &   13.36       \\ 			
		$R_{k}$     &   1.58e-11    & 4.50e-11  &  2.96e-11   &  2.58e-11  &  1.05e-10   \\ 
		$E_k$&  4.66e-11     & 1.11e-10 &  5.24e-11   &  4.30e-11   &  1.10e-10 \\ \hline
		
		$\alpha$& 6& 7  &8 & 9 & 10  \\
		$n=65536$  & & & & & \\ \hline
		
		Iter &   13     &   14    &   14   &  14    & 15   \\
		CPU &   8.96  &  8.88    &  8.91  &  9.38  &  10.28        \\ 			
		$R_{k}$     &   8.57e-11    &  9.00e-12  &  1.88e-11    &  5.54e-11  &  3.69e-12    \\ 
		$E_k$&    1.03e-10   & 4.26e-11 &  7.30e-11   &  5.28e-11   & 2.76e-10  \\ \hline
		
	\end{tabular}
\end{table}

\begin{table}[!ht]
	
	\caption{Average number of the iterations of Chebyshev semi-iteration for solving the inner systems with the PRESB preconditioner versus the inner tolerance for Example \ref{Ex1} and preconditioner I with $\omega=1$ and $m=256$. \label{Tbl11}}
	
	\begin{tabular}{|c|c|c|} \hline
		inner tolerance & $W_1+iT$  &  $W_2-iT$  \\ \hline
		$10^{-2}$       &    2      &  4             \\ 
		$10^{-3}$       &    4      &  5.43           \\ 
		$10^{-4}$       &    5      &  6.85           \\ 
		$10^{-5}$       &    6      &  8.14           \\ 
		$10^{-6}$       &    8      &  9.42           \\ 
		$10^{-7}$       &    9      &  10.85           \\ 
		$10^{-8}$       &    10     &  12.14           \\ 
		$10^{-9}$       &    12     &  13.85           \\ 
		$10^{-10}$      &    13     &  15.14          \\ \hline 
	\end{tabular}
\end{table}

\section{Conclusion}\label{sec04}

We have presented three iteration methods for solving a class of  complex symmetric system of linear equations $(W+iT)x=b$, where $W$ is indefinite and $T$ is symmetric positive definite. We have proved that they converge unconditionally. The induced preconditioners  have been applied to accelerate the convergence of the GMRES method for solving the system. Efficient ways to implement the preconditioners have been presented.  Numerical experiments show that the preconditioners can compete with the recently presented method, SNSS.

\end{document}